# DETECTION OF SPATIAL CLUSTERING WITH AVERAGE LIKELIHOOD RATIO TEST STATISTICS

By Hock Peng Chan[1]

*National University of Singapore*

Generalized likelihood ratio (GLR) test statistics are often used in the detection of spatial clustering in case-control and case-population datasets to check for a significantly large proportion of cases within some scanning window. The traditional spatial scan test statistic takes the supremum GLR value over all windows, whereas the average likelihood ratio (ALR) test statistic that we consider here takes an average of the GLR values. Numerical experiments in the literature and in this paper show that the ALR test statistic has more power compared to the spatial scan statistic. We develop in this paper accurate tail probability approximations of the ALR test statistic that allow us to by-pass computer intensive Monte Carlo procedures to estimate *p*-values. In models that adjust for covariates, these Monte Carlo evaluations require an initial fitting of parameters that can result in very biased *p*-value estimates.

**1. Introduction.** The detection of local clustering in spatial point processes is of interest in epidemiological studies, forestry, geological studies, neural imaging and astronomy. There are a number of excellent texts and review papers on this, including [5, 13, 29]. A classical application that will be used here as an illustrative example is the identification of potential sources of environmental pollution that have contributed to higher rates of disease cases for residents living in their vicinity.

Let $\mathcal{T} = \{\mathbf{t}_i : 1 \leq i \leq I\}$, with $\mathbf{t}_i \in \mathbf{R}^d$ denoting the location of the $i$th case. We are interested in the presence of an unusually large number of cases near an unspecified location $\mathbf{v} = (v_1, \ldots, v_d)$ inside a bounded domain $D$. If $\mathcal{T}$ is generated from a process with known and constant intensity under the null hypothesis, we can test for the presence of clusters by computing the

Received November 2008; revised March 2009.
[1]Supported by grants from the National University of Singapore.

maximal number of cases in the cubic windows $\prod_{k=1}^{d}[v_k - \frac{w}{2}, v_k + \frac{w}{2}]$, over all $\mathbf{v} \in D$ for a fixed window size $w > 0$. The question of whether this number is significantly large or may have occurred with reasonable chance under the null hypothesis was addressed in [21, 23], via asymptotic $p$-value calculations and $p$-value bounds. Extensions to weighted counting using kernel functions were also achieved in [27].

Rather than assuming that the underlying intensity is known and constant, we can assume instead that a control dataset $\mathcal{U} = \{\mathbf{u}_j : 1 \leq j \leq J - I\}$ is available for estimation of the possibly nonconstant intensity function. There has been considerable work done on the use of kernel functions to smooth $\mathcal{U}$ to provide an intensity estimate, and the significance of a cluster of cases is calculated by assuming that the estimated intensity is the true intensity (see, e.g., [1, 7, 9] and references therein). An alternative approach, as considered in [6, 26], is to merge $\mathcal{T}$ and $\mathcal{U}$ into a combined dataset $\mathcal{X} := \{(\mathbf{t}_i, 1) : 1 \leq i \leq I\} \cup \{(\mathbf{u}_j, 0) : 1 \leq j \leq J - I\}$ and rewrite it as $\{(\mathbf{x}_i, X_i) : 1 \leq i \leq J\}$. The SaTScan software developed by Kulldorff and Information Management Services Inc. [16] (see also [17]) considers merged datasets, with generalized likelihood ratio (GLR) test statistics used to provide a score for each window, and the spatial scan statistic, the supremum GLR score used to determine significance. Instead of cubic windows, spherical windows $C(\mathbf{v}, w) := \{\mathbf{t} : \sum_{k=1}^{d}(v_k - t_k)^2 \leq w^2\}$ are considered.

In Section 2, we consider the average likelihood ratio (ALR) test statistic, which uses an average rather than the supremum GLR score as the summary test statistic. Numerical studies in the literature and in this paper show that the ALR test statistic has more power compared to the spatial scan test statistic. We provide moderate deviation tail probability approximations in Section 2.1 for the ALR test statistic and illustrate their extensions to logistic regression models for covariate adjustments in Section 3. These $p$-value approximations allow us to avoid the use of computationally expensive Monte Carlo methods and are especially important when covariate adjustments are required, as the Monte Carlo method currently in use requires an initial fitting of parameters that can result in very biased $p$-value estimates (see Examples 1 and 2 in Section 3.1). In Section 4, we perform comparison studies on real and simulated datasets. A discussion is provided in Section 5 followed by derivations of the asymptotic formulae in Section 6. The appendices contain technical details and proofs.

**2. The spatial scan and ALR test statistics.** Throughout this paper, we shall use $\| \cdot \|$ to denote the $L_2$ norm of a vector. For any set $A$, vector $\mathbf{t}$ and real number $b$, we shall let $\mathbf{t} + bA = \{\mathbf{t} + b\mathbf{a} : \mathbf{a} \in A\}$. We shall use $\mathbf{I}$ to denote the indicator function and $\#$ to denote the number of elements in a finite set. For constants $a_n$ and $b_n$, the notation $a_n \sim b_n$ shall mean $a_n/b_n \to 1$, while for random variables $Y_1, Y_2, \ldots$ and $Z_1, Z_2, \ldots$, the notation $Y_n \sim Z_n$ shall mean $Y_n/Z_n \xrightarrow{p} 1$. We shall use $\mathbf{Z}$ to denote the set of integers and $\mathbf{0}$ to denote the zero vector. We shall also adopt the conventions $0 \log 0 = 0$ and $0^0 = 1$.

Let $\mathcal{X} = \{(\mathbf{x}_i, X_i) : 1 \leq i \leq J\}$, where $\mathbf{x}_i$ denotes the location of the $i$th subject, while $X_i = 1$ if the subject is a case and $X_i = 0$ otherwise. Conditioned on $\mathbf{x} := (\mathbf{x}_1, \ldots, \mathbf{x}_J)$, the random vector $\mathbf{X} := (X_1, \ldots, X_J)$ consists of independent Bernoulli random variables. Under the null hypothesis $H_0$ of no clustering, there exists $p_0 \in (0, 1)$ such that

$$P_0\{X_i = 1\} = p_0 \qquad \text{for all } 1 \leq i \leq J. \tag{2.1}$$

Let $B$ be a subset of $\mathbf{R}^d$ and $H_B^{(1)}$ the hypothesis that there exists $p_1 > p_2$ such that

$$\begin{aligned} &P\{X_i = 1 | \mathbf{x}_i \in B\} = p_1, \\ &P\{X_i = 1 | \mathbf{x}_i \notin B\} = p_2 \qquad \text{for all } 1 \leq i \leq J. \end{aligned} \tag{2.2}$$

Let $\widehat{p}_0 = I/J$ be the maximum likelihood estimate (MLE) of $p_0$ under $H_0$ and let

$$\phi(p) = p \log\left(\frac{p}{\widehat{p}_0}\right) + (1 - p) \log\left(\frac{1 - p}{1 - \widehat{p}_0}\right). \tag{2.3}$$

Let $m_B = \sum_{i=1}^{J} \mathbf{I}_{\{\mathbf{x}_i \in B, X_i = 1\}}$ be the number of cases and $n_B = \sum_{i=1}^{J} \mathbf{I}_{\{\mathbf{x}_i \in B\}}$ the number of subjects in $B$. The log GLR score for testing $H_0$ against $H_B^{(1)}$ is

$$\begin{aligned} S^{(1)}(B) &:= \log\Big\{ \sup_{1 \geq p_1 > p_2 \geq 0} [p_1^{m_B}(1 - p_1)^{n_B - m_B} p_2^{I - m_B}(1 - p_2)^{J - I - (n_B - m_B)}] \Big\} \\ &\quad - \log[\widehat{p}_0^I (1 - \widehat{p}_0)^{J - I}] \\ &= \left[ n_B \phi\left(\frac{m_B}{n_B}\right) + (J - n_B)\phi\left(\frac{I - m_B}{J - n_B}\right) \right] \mathbf{I}_{\{m_B/n_B > \widehat{p}_0\}}. \end{aligned}$$

To detect both over- and under-clustering, we compare $H_0$ against the two-sided alternative hypothesis $H_B^{(2)}$ that (2.2) holds for some $p_1 \neq p_2$. The log GLR score is then

$$S^{(2)}(B) := n_B \phi\left(\frac{m_B}{n_B}\right) + (J - n_B)\phi\left(\frac{I - m_B}{J - n_B}\right). \tag{2.4}$$

Let $\mathcal{B}$ be a finite class of measurable subsets of $\mathbf{R}^d$, possibly dependent on $\mathbf{x}$ but not on $\mathbf{X}$. The spatial scan statistic for testing $H_0$ vs. $\bigcup_{B \in \mathcal{B}} H_B^{(k)}$, $k = 1$ or $2$, is

$$M_{\mathcal{B}}^{(k)} := \sup_{B \in \mathcal{B}} S^{(k)}(B). \tag{2.5}$$

The spatial scan statistic has the drawback of not making full use of information provided by secondary clusters to conclude the presence of local clustering. For example, if there are scores $S^{(k)}(B_1) > S^{(k)}(B_2)$ for nonoverlapping windows $B_1$ and $B_2$ both slightly smaller than the critical value, the information provided by $S^{(k)}(B_2)$ is not utilized in the decision not to reject $H_0$. Gangnon and Clayton [12] introduced the weighted ALR test statistic

$$\sum_{B \in \mathcal{B}} w_B e^{S^{(2)}(B)} \qquad \text{with } w_B > 0 \text{ for all } B \in \mathcal{B} \quad \text{and} \quad \sum_{B \in \mathcal{B}} w_B = 1.$$

Unlike the spatial scan statistic, significance for the weighted ALR test statistic can be concluded based on many moderately large scores. The numerical studies in [12] suggest that the weighted ALR is more powerful than the spatial scan statistic in the detection of local clusters. Siegmund [28] also reports a closely related test statistic that is slightly more powerful, compared to the scan test statistic in a numerical study on the genome scan. This is in contrast to global clustering test statistics like $(\#\mathcal{B})^{-1} \sum_{B \in \mathcal{B}} S^{(2)}(B)$, which are expected to have lower power compared to the spatial scan statistic when only a few local clusters are present (see [18] for supporting numerical results). We consider in this paper $p$-value approximations for the (log) ALR test statistic

$$(2.6) \qquad U_{\mathcal{B}}^{(k)} := 2 \log \left( (\#\mathcal{B})^{-1} \sum_{B \in \mathcal{B}} e^{S^{(k)}(B)} \right).$$

An extension of these approximations to weighted ALR test statistics is given in the appendices of [4].

2.1. *Moderate deviation tail probabilities.* In this paper, we provide tail approximations of the ALR test statistics under the following assumptions.

(A1) The domain $D$ is a compact subset of $\mathbf{R}^d$ and satisfies

$$\#\{\mathbf{t} \in (\varepsilon \mathbf{Z})^d : t + [0, \varepsilon]^d \subset D\} \sim \#\{\mathbf{t} \in (\varepsilon \mathbf{Z})^d : (\mathbf{t} + [0, \varepsilon]^d) \cap D \neq \varnothing\} \sim |D|/\varepsilon^d$$

as $\varepsilon \to 0$.

(A2) The locations $\mathbf{x}_1, \ldots, \mathbf{x}_J$ are independent and identically distributed (i.i.d.) random vectors generated from $\lambda$, a continuous and positive density on $D$.

(A3) The class of scanning sets $\mathcal{B}$ is a sub-class of $\mathcal{C} := \{\mathbf{v} + wA : \mathbf{v} \in D, w_0 \leq w \leq w_1\}$, where $A$ is a convex, open and bounded subset of $\mathbf{R}^d$, with $\mathbf{0} \in A$ and $0 < w_0 \leq w_1 < (|D|/|A|)^{1/d}$.

In Theorem 1 below and Theorem 2 in Section 3, $\mathcal{B}(=\mathcal{B}_c)$ may vary with the critical value $c$ and constraints are placed only on the growth of $J$ (for Theorem 1) and $\#\mathcal{B}$ with respect to $c$. The class of $\mathcal{C}$ of candidate scanning sets is, however, fixed for all $c > 0$. The proofs of the theorems use change of measure arguments and linearization techniques developed by Lai and Siegmund [19, 20] and Woodroofe [32, 33], to analyze GLR test statistics in sequential analysis and are given in Section 5. A motivation of the proofs is also given by a simpler Theorem 3 and its proof in Appendix A. Let $\chi_1^2$ denote a chi-square random variable with one degree of freedom.

THEOREM 1.  *Assume* (A1)–(A3) *and let* (2.1) *hold for some* $0 < p_0 < 1$. *Let* $\log(\#\mathcal{B}) = o(c^{1/3})$ *and assume that* $c \sim \kappa J^s$ *for some* $\kappa > 0$ *and* $0 < s < 1$. *Then as* $c \to \infty$,

$$(2.7) \qquad P_0\{U_{\mathcal{B}}^{(k)} \geq c | \mathbf{x}\} \sim k P\{\chi_1^2 \geq c\}/2 \qquad \text{for } k = 1, 2.$$

The assumptions (A2), (A3) and the relation $c \sim \kappa J^s$ in the statement of Theorem 1 are needed to ensure that the number of subjects in each $B \in \mathcal{B}$ approaches infinity fast enough for a chi-square tail probability approximation of $S^{(k)}(B)$ to hold. This leads to the chi-square tail probability approximation of $U_{\mathcal{B}}^{(k)}$. The uniform approximation when conditioning on $\mathbf{x}$ in (2.7) ensures that we do not reject $H_0$ unevenly with respect to the configuration of the locations. However, it is also important for us to check the actual type I error probability when $\mathbf{x}$ is not conditioned on (see Example 2 in Section 3.1).

**3. Logistic modeling.** To see why (2.7) extends to more complicated models, it is useful to view it as resulting from two different asymptotics. Let $\lambda_B = \int_B \lambda(\mathbf{t}) \, d\mathbf{t}$, where $\lambda$ is the density in (A2). Let $\omega$ be Gaussian white noise with $\omega(B) \sim N(0, \lambda_B)$ for $B \subset D$ and $\omega(A), \omega(B)$ independent whenever $A$ and $B$ are disjoint. Let $Z_B = \lambda_B^{-1/2}(1 - \lambda_B)^{-1/2}[\omega(B) - \lambda_B\omega(D)]$. The first asymptotic is a weak convergence of $S^{(2)}(B)$ to $Z_B^2/2$ uniformly over $B \in \mathcal{C}$, and this holds largely because $\inf_{B \in \mathcal{C}}(n_B/c) \to \infty$ when $c \sim \kappa J^s$ for $0 < s < 1$. The second asymptotic is like (2.7) [see (3.2) below], but with ALRs $U_{\mathcal{B}}^{(2)}$ and $U_{\mathcal{B}}^{(1)}$ replaced by

$$(3.1) \qquad \begin{aligned} U_Z^{(2)} &:= 2\log\left((\#\mathcal{B})^{-1}\sum_{B \in \mathcal{B}} e^{Z_B^2/2}\right) \quad \text{and} \\ U_Z^{(1)} &:= 2\log\left((\#\mathcal{B})^{-1}\sum_{B \in \mathcal{B}} e^{Z_{B+}^2/2}\right), \end{aligned}$$

respectively, where $Z_{B+} = \max\{Z_B, 0\}$.

THEOREM 2.  *Assume* (A1), (A3) *and let* $\log(\#\mathcal{B}) = o(c^{1/3})$. *Then as* $c \to \infty$,

$$(3.2) \qquad P\{U_Z^{(k)} \geq c\} \sim k P\{\chi_1^2 \geq c\}/2 \qquad \text{for } k = 1, 2.$$

Since $|U_Z^{(2)} - U_{\mathcal{B}}^{(2)}| \leq 2\sup_{B \in \mathcal{B}} |S^{(2)}(B) - Z_B^2/2|$ and $|U_Z^{(1)} - U_{\mathcal{B}}^{(1)}| \leq 2 \times \sup_{B \in \mathcal{B}} |S^{(1)}(B) - Z_{B+}^2/2|$, the combination of the two asymptotics described above provides us with chi-square tail approximations for $U_{\mathcal{B}}^{(k)}$.

Consider more generally datasets containing additional information like the age, sex, diet and smoking habits of the subjects. These covariates may influence the outcome and, hence, we may have to correct for spatial imbalances of these covariates when testing for spatial clustering. Let $\mathbf{u}_i = (u_{i1}, \ldots, u_{ir})'$ be the covariate vector of the $i$th subject, with $u_{i1} = 1$ denoting the intercept term, and let $p_i = P\{X_i = 1|\mathbf{x}_i, \mathbf{u}_i\}$. Consider the logistic model

$$(3.3) \qquad\qquad p_i = (1 + e^{-\beta'\mathbf{u}_i - \theta_i})^{-1},$$

where $\beta = (\beta_1, \ldots, \beta_r)'$ is a nuisance parameter vector. Under the null hypothesis $H_0$ of no clustering, $\theta_i = 0$ for all $i$, while under the one-sided alternative hypothesis $H_B^{(1)}$, $\theta_i = \theta\mathbf{I}_{\{\mathbf{x}_i \in B\}}$ for some $\theta > 0$. Under the two-sided alternative hypothesis $H_B^{(2)}$, $\theta_i = \theta\mathbf{I}_{\{\mathbf{x}_i \in B\}}$ for some $\theta \neq 0$. Let $\widehat{\beta}$ be the MLE of $\beta$ under $H_0$ and $(\widehat{\beta}_B^{(k)}, \widehat{\theta}_B^{(k)})$ the MLE of $(\beta, \theta)$ under $H_0 \cup H_B^{(k)}$. Define

$$(3.4) \qquad \begin{aligned} &\widehat{p}_i = (1 + e^{-\widehat{\beta}'\mathbf{u}_i})^{-1}, \qquad \widehat{p}_{iB}^{(k)} = (1 + e^{-\widehat{\beta}_B^{(k)'}\mathbf{u}_i - \widehat{\theta}_B^{(k)}\mathbf{I}(\mathbf{x}_i \in B)})^{-1}, \\ &Y_{iB}^{(k)} = X_i \log\left(\frac{\widehat{p}_{iB}^{(k)}}{\widehat{p}_i}\right) + (1 - X_i)\log\left(\frac{1 - \widehat{p}_{iB}^{(k)}}{1 - \widehat{p}_i}\right). \end{aligned}$$

Then the ALR test statistics are

$$(3.5) \quad U_{\mathcal{B}}^{(k)} = 2\log\left((\#\mathcal{B})^{-1}\sum_{B \in \mathcal{B}} e^{S^{(k)}(B)}\right) \qquad \text{where } S^{(k)}(B) = \sum_{i=1}^{J} Y_{iB}^{(k)}.$$

The scores $S^{(k)}(B)$ are asymptotically chi-square, even when $\beta$ is infinite dimensional (see [2, 22] and references therein). The efficient score expansions of the log profile likelihoods that are used for deriving these chi-square approximations can also be used to provide the covariance structure of the limiting multivariate normal of $\sqrt{n}\widehat{\theta}_B^{(2)}$ over $B \in \mathcal{B}$, and this structure depends on the nuisance parameter under $H_0$ (see Appendix B for more details). However, the chi-square approximations of $U_{\mathcal{B}}^{(k)}$ in the moderate deviations domain do not depend on the covariance structure of the limiting multivariate normal. In other words,

$$(3.6) \qquad P_{(0,\beta)}\{U_{\mathcal{B}}^{(k)} \geq c\} \sim kP\{\chi_1^2 \geq c\}/2 \qquad \text{for } k = 1, 2$$

uniformly over compact sets of $\beta$ (see Appendix B).

This is desirable because the $p$-value is in principle computed from the worst-case scenario under $H_0$. In this respect, the ALR test statistic shares the same uniform asymptotics as the GLR test statistic for a composite null hypothesis versus single composite alternative hypothesis with a dimension difference of one, differing only in that for the GLR test statistic, the approximation occurs in the central limit domain as well. The spatial scan test statistic does not have such uniform asymptotics over nuisance parameters. Hence Theorems 1 and 2 are not just devices for $p$-value approximations, but also theoretical results that provide understanding of the asymptotic properties of the ALR test statistic. To reduce computational time for large datasets, we can avoid searching for a new $(\widehat{\beta}_B^{(k)}, \widehat{\theta}_B^{(k)})$ for each $B \in \mathcal{B}$ by replacing $S^{(2)}(B)$ by a first-order quadratic approximation (see either (4)–(6) of [22] or (B.1) in Appendix B).

3.1. *Monte Carlo evaluation of conditional $p$-values.* Under (2.1), the conditional $p$-value $P_0\{M_{\mathcal{B}}^{(k)} \geq c | I, \mathbf{x}\}$ does not depend on $p_0$ and can be evaluated by a permutation test. Permutation tests are nonparametric tests that compute $p$-values from permutations of the observations $X_1, \ldots, X_J$, which are often assumed to be i.i.d. under the null hypothesis. In principle, the $p$-value is the fraction of permutations with values of test statistics at least as large as the original test statistic, though in practice the number of permutations is usually too large for direct computations, and Monte Carlo methods are used instead to sample a random subset of permutations for $p$-value estimation (for more details, see [10, 11]). In the SaTScan software, users are prompted to select $L = 99$, 999 or 9999 random permutations. For each $1 \leq \ell \leq L$, compute $M_{\mathcal{B}, \ell}^{(k)}$ from $\{(\mathbf{x}_i, X_{i\ell}) : 1 \leq i \leq J\}$, where $(X_{1\ell}, \ldots, X_{J\ell})$ is a random permutation of $(X_1, \ldots, X_J)$. Then the estimated conditional $p$-value is $(1 + \sum_{\ell=1}^{L} \mathbf{I}_{\{M_{\mathcal{B}, \ell}^{(k)} \geq c\}})/(1 + L)$. The extension of the method to estimate $P_0\{U_{\mathcal{B}}^{(k)} \geq c | I, \mathbf{x}\}$ is straightforward.

When covariates are present, the SaTScan software uses the following Monte Carlo procedure, as advocated in [15]. Assume that there are $n_j$ subjects at location $\mathbf{v}_j$ for $1 \leq j \leq q$, with $n_j$ large. Fit (3.3) under the null hypothesis $H_0$, that there are no spatial effects, that is with $\theta_i = 0$ for all $i$. The fitted value $\widehat{p}_i$, given in (3.4), is the estimated risk of the $i$th subject. At each $\mathbf{v}_j$, estimate the total risk by $\eta_j = \sum_{i : \mathbf{x}_i = \mathbf{v}_j} \widehat{p}_i$. Let $m_j = \sum_{i : \mathbf{x}_i = \mathbf{v}_j} X_i$, $m_B = \sum_{\mathbf{v}_j \in B} m_j$ and $\eta_B = \sum_{\mathbf{v}_j \in B} \eta_j$. Assume that under $H_0$, $m_1, \ldots, m_q$ are independent Poisson random variables with respective means $\eta_1, \ldots, \eta_q$. Then conditioned on $m_1 + \cdots + m_q = I$, the adjusted spatial scan statistic for testing $H_0$ against $\bigcup_{B \in \mathcal{B}} H_B^{(2)}$ is

$$\widetilde{M}_{\mathcal{B}}^{(2)} := \sup_{B \in \mathcal{B}} \widetilde{S}^{(2)}(B) \qquad \text{where}$$

TABLE 1
*Comparison of the type I error probabilities and detection powers of $\widetilde{M}_{\mathcal{B}}^{(2)}$ and $U_{\mathcal{B}}^{(2)}$ at significance levels $\alpha = 0.05$ and $\alpha = 0.01$ with 1000 independent copies of $\mathbf{X}$*

| | $\alpha = 0.05$ | | $\alpha = 0.01$ | |
|---|---|---|---|---|
| $\theta$ | MC: $\widetilde{M}_{\mathcal{B}}^{(2)}$ | ALR: $U_{\mathcal{B}}^{(2)}$ | MC: $\widetilde{M}_{\mathcal{B}}^{(2)}$ | ALR: $U_{\mathcal{B}}^{(2)}$ |
| 0 | 0.026 | 0.048 | 0.004 | 0.008 |
| 0.2 | 0.088 | 0.158 | 0.021 | 0.054 |
| 0.4 | 0.367 | 0.499 | 0.137 | 0.261 |
| 0.6 | 0.740 | 0.849 | 0.506 | 0.676 |

$$(3.7) \qquad \widetilde{S}^{(2)}(B) := m_B \log\left(\frac{m_B}{\eta_B}\right) + (I - m_B) \log\left(\frac{I - m_B}{I - \eta_B}\right).$$

To simulate the Monte Carlo $p$-value for each $1 \leq \ell \leq L$, where $L$ is the required number of simulation runs, generate $(m_{1\ell}, \ldots, m_{q\ell})$ from a multinomial distribution with $I$ trials and success probabilities $(\eta_1/I, \ldots, \eta_q/I)$, then compute $\widetilde{M}_{\mathcal{B},\ell}^{(2)}$ using (3.7). The estimated $p$-value is then

$$\left(1 + \sum_{\ell=1}^{L} \mathbf{I}_{\{\tilde{M}_{\mathcal{B},\ell}^{(2)} \geq \tilde{M}_{\mathcal{B}}^{(2)}\}}\right) \Big/ (1 + L).$$

EXAMPLE 1. Let $D$ be a union of disjoint sets $B_1$, $B_2$ and $B_3$, each containing 1000 subjects. Generate dummy covariates $u_i \sim N(0,1)$ if $\mathbf{x}_i \in B_2 \cup B_3$ and $u_i \sim N(1,1)$ if $\mathbf{x}_i \in B_1$, then keep them fixed for the remaining part of this exercise. Let $\mathcal{B} = \{B_1, B_2, B_3\}$ and let

$$(3.8) \qquad P\{X_i = 1 | \mathbf{x}_i, u_i\} = (1 + e^{-\beta_1 - \theta \mathbf{I}_{\{\mathbf{x}_i \in B_1\}}})^{-1}.$$

In our comparison study, we generate $\mathbf{X} = (X_1, \ldots, X_{3000})$ from (3.8) with $\beta_1 = -3$, $\theta \geq 0$ and compute the Monte Carlo $p$-values of $\widetilde{M}_{\mathcal{B}}^{(2)}$ with $L = 999$ simulation runs and also the $p$-values of $U_{\mathcal{B}}^{(2)}$ using chi-square tail probability approximations. The scores $S^{(2)}(B)$ are computed from (3.4)–(3.5) with $u_i$, the only covariate of the $i$th subject. For each $\theta \geq 0$, the above procedure is repeated 1000 times, each time with a different copy of $\mathbf{X}$. The estimated type I error probabilities and power are summarized in Table 1. We see that the Monte Carlo risk adjustment method provides very conservative $p$-values (see [3] for alternative strategies to deal with this drawback).

EXAMPLE 2. We choose a slightly different design here to check the type I error probability and power $P\{U_{\mathcal{B}}^{(2)} \geq c\}$ (without conditioning on $\mathbf{x}$). In each simulation run, twenty locations $\mathbf{v}_1, \ldots, \mathbf{v}_{20}$ are generated uniformly and randomly on the unit square $[0,1]^2$. Let $C_0$ be a circle of radius 0.3, centered at $(0.5, 0.5)$. Fifty subjects are located at each $\mathbf{v}_i$, each of them generated as a case with probability $p_0 = 0.05$ if $\mathbf{v}_i \notin C_0$, and generated as a case with probability $p_1 \geq p_0$ if $\mathbf{v}_i \in C_0$. Each subject at $\mathbf{v}_i$ is given a dummy covariate distributed as $N(0,1)$ if $\mathbf{v}_i \notin C_0$ and distributed as $N(1,1)$ if $\mathbf{v}_i \in C_0$. For each $1 \leq i \leq 20$, let $0 = r_{i,1} < \cdots < r_{i,20}$ be the ordered values of $\|\mathbf{v}_j - \mathbf{v}_i\|$ for $j = 1, \ldots, 20$. We consider the class of scanning sets

$$\mathcal{B} = \{C(\mathbf{v}_i, r_{i,j}) : 1 \leq i \leq 20, 1 \leq j \leq 10\},$$

where $C(\mathbf{v}, r)$ is a circle of radius $r$, centered at $\mathbf{v}$. One thousand simulation runs are used to estimate each type I error probability and power of the adjusted scan statistic $\widetilde{M}_{\mathcal{B}}^{(2)}$ (using $L = 999$ permutations) and the ALR test statistic $U_{\mathcal{B}}^{(2)}$ (using the chi-square distribution) (see Table 2). We see that the Monte Carlo method has low type I error probability and corresponding loss of power when compared against the ALR test statistic.

TABLE 2
*Comparison of the type I error probabilities and detection powers of $\widetilde{M}_{\mathcal{B}}^{(2)}$ and $U_{\mathcal{B}}^{(2)}$ at significance levels $\alpha = 0.05$ and $\alpha = 0.01$ with 1000 simulation runs*

| $p_1$ | $\alpha = 0.05$ | | $\alpha = 0.01$ | |
|---|---|---|---|---|
| | MC: $\widetilde{M}_{\mathcal{B}}^{(2)}$ | ALR: $U_{\mathcal{B}}^{(2)}$ | MC: $\widetilde{M}_{\mathcal{B}}^{(2)}$ | ALR: $U_{\mathcal{B}}^{(2)}$ |
| 0.05 | 0.026 | 0.053 | 0 | 0.007 |
| 0.2 | 0.054 | 0.119 | 0.008 | 0.026 |
| 0.4 | 0.243 | 0.395 | 0.111 | 0.209 |
| 0.6 | 0.514 | 0.682 | 0.327 | 0.484 |

**4. Numerical examples.** We analyze a case-control dataset in Section 4.1, a case-population dataset in Section 4.2 and various simulated datasets in Section 4.3.

4.1. *Laryngeal cancer dataset.* This dataset consists of: (i) the locations of 58 cases of laryngeal cancer occurring in two districts in Lancashire for the period 1974–1985; and (ii) the locations of 978 control cases of lung cancer for the same period and districts in the domain $D = [34500, 36500] \times [41100, 43100]$ (see [8] for more background). A key feature is a cluster of four laryngeal cancer cases (see the bottom of the left plot of Figure 1) located near an industrial waste incinerator, which is considered a potential source of the cluster of laryngeal cancer cases. We want to test for the presence of local clusters without biasing ourselves a priori with information on the possible sources of the laryngeal cancer cases. As the location co-ordinates in the datasets are rounded to the nearest tens, we consider the covering sets

$$\mathcal{B}_w := \{C(\mathbf{v}, w) : \mathbf{v} \in D \cap (10\mathbf{Z} + 5)^2, n_{C(\mathbf{v},w)} \geq 2\}$$

with radii $w = 40$, 50, 60 and 70. Hence each circle in $\mathcal{B}_w$ contains at least two subjects and has a center $\mathbf{v}$ with co-ordinates ending with 5 and lies inside $D$. Express $M^{(1)}_{\mathcal{B}_w}$ and $U^{(1)}_{\mathcal{B}_w}$ more simply as $M^{(1)}_w$ and $U^{(1)}_w$, respectively.

In Table 3, we tabulate Monte Carlo conditional $p$-values of both $M^{(1)}_w$ and $U^{(1)}_w$ using the permutation method described in Section 3.1. We observe that for both the spatial scan and ALR test statistics, $p$-values below 0.02 are obtained when $w = 40$. This is in contrast to $p$-values of 0.08 to 0.8 obtained using kernel-based methods (see [1]). The choice of window size $w$ affects the $p$-value substantially when using $M^{(1)}_w$, and this is also true when using kernel-based methods. In contrast, the influence of window size on the $p$-values of $U^{(1)}_w$ is much smaller. In this sense, the ALR test statistic is more robust against misspecification of cluster shape and size, that is, when $H^{(1)}_B$ is true for some $B \notin \mathcal{B}$, because under such a situation there will often be many windows having moderately large scores, and this will aid the rejection of $H_0$. The construction of Table 3 requires a substantial amount of computation as there are more than 5000 scanning sets in each $\mathcal{B}_w$.

A numerical power study (see Table 4) indicates that the ALR and spatial scan test statistics do not dominate each other when there is only one source of spatial clustering. In this study, we fix the locations $\mathbf{x}$ and the total number of cases $I = 58$. Consider a circle with radius 40 and let $n$ be the number of points in it. Let $p$ be the probability that a point in the circle is

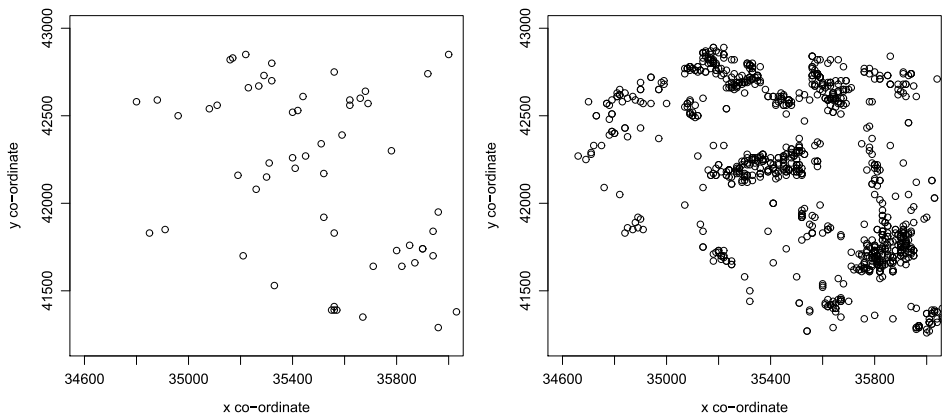

FIG. 1. *Scatter plots of the 58 laryngeal cancer cases (left) and the 978 lung cancer cases (right).*

*Numerical values of the test statistics and Monte Carlo conditional p-value estimates $\pm$ standard error for $B_w$ with 2000 simulation runs for each spatial scan statistic p-value and 10,000 runs for each ALR test statistic p-value*

| | Spatial scan statistic $M_w^{(1)}$ | | ALR test statistic $U_w^{(1)}$ | |
|---|---|---|---|---|
| $w$ | Value | MC p-val. (cond.) | Value | MC p-val. (cond.) |
| 40 | 9.21 | $0.016 \pm 0.003$ | 5.29 | $0.0104 \pm 0.0010$ |
| 50 | 7.95 | $0.090 \pm 0.006$ | 4.47 | $0.0137 \pm 0.0012$ |
| 60 | 7.95 | $0.078 \pm 0.006$ | 4.07 | $0.0200 \pm 0.0014$ |
| 70 | 7.95 | $0.079 \pm 0.006$ | 3.89 | $0.0213 \pm 0.0014$ |

simulated as a case and $\widetilde{p}$ the probability that a point outside the circle is simulated as a case. Thus the relative ratio (RR) is $p/\widetilde{p}$. The numbers $p$ and $\widetilde{p}$ are determined from the constraint

$$np + (1036 - n)\widetilde{p} = 58.$$

In the $\ell$th simulation run, $1 \le \ell \le 1000$, we generate $\{X_{i\ell} : 1 \le i \le 1036\}$ with success probabilities $p$ (for $\mathbf{x}_i$ inside circle) or $\widetilde{p}$ (for $\mathbf{x}_i$ outside circle), and repeat until a total of 58 cases is observed before proceeding to compute $U_{40,\ell}^{(1)}$ and $M_{40,\ell}^{(1)}$. The estimated power is the proportion of runs in which the critical value is equaled or exceeded.

We also try out scanning sets with different radii at different centers as suggested by a referee, and obtain similar $p$-values for the spatial scan and ALR test statistics (see Table 5). The classes of scanning sets considered here are of the form

$$\mathcal{B}_j = \{C(\mathbf{x}_i, r_{ij}) : 1 \le i \le 1036\} \qquad \text{for } j = 5, 6, 7,$$

where $0 = r_{i1} \le r_{i2} \le \cdots$ are the ordered values of $\|\mathbf{x}_i - \mathbf{x}_k\|$ for $1 \le k \le 1036$. It is interesting to note that even though the largest window score of 9.21,

*Powers of $U_{40}^{(1)}$ and $M_{40}^{(1)}$ based on 1000 simulation runs for each entry, with estimated 1% critical values $c_{M,0.01} = 9.49$ and $c_{U,0.01} = 5.29$. In each row, the n points lying in a circle centered at $(v_1, v_2)$ with radius $w = 40$ are simulated as cases with probability RR times larger than points lying outside the circle*

| $v_1$ | $v_2$ | $n$ | RR | Power of $U_{40}^{(1)}$ | Power of $M_{40}^{(1)}$ |
|---|---|---|---|---|---|
| 35565 | 41395 | 6 | 12 | $0.49 \pm 0.02$ | $0.36 \pm 0.02$ |
| 35195 | 42745 | 9 | 11 | $0.54 \pm 0.02$ | $0.54 \pm 0.02$ |
| 35515 | 42255 | 12 | 10 | $0.52 \pm 0.02$ | $0.68 \pm 0.01$ |
| 35255 | 42155 | 15 | 8 | $0.45 \pm 0.02$ | $0.47 \pm 0.02$ |
| 35595 | 42745 | 18 | 7 | $0.47 \pm 0.02$ | $0.42 \pm 0.02$ |

*Numerical values of the test statistics and Monte Carlo conditional p-value estimates $\pm$ standard error with $2000$ simulation runs in each entry for scanning sets with different radii*

| $j$ | Spatial scan statistic $M^{(1)}$ | | ALR test statistic $U^{(1)}$ | |
|---|---|---|---|---|
| | **Value** | **MC p-val. (cond.)** | **Value** | **MC p-val. (cond.)** |
| 5 | 9.21 | $0.016 \pm 0.003$ | 5.38 | $0.012 \pm 0.002$ |
| 6 | 7.95 | $0.043 \pm 0.005$ | 5.76 | $0.006 \pm 0.002$ |
| 7 | 7.04 | $0.079 \pm 0.006$ | 3.70 | $0.027 \pm 0.004$ |

obtained from a scanning set containing four cases and one control, is missed when $j = 6$, the ALR score actually increased.

4.2. *New York leukaemia dataset.* We use here an updated version of the dataset presented in [31], which tracks leukaemia occurrences in 281 census tracts in New York state. Let $\mathbf{v}_j$ denote the centroid of the $j$th census tract and let $m_j$ and $n_j$ be the number of leukaemia cases and population size, respectively, at $\mathbf{v}_j$. Let $m_B = \sum_{\mathbf{v}_j \in B} m_j$, $n_B = \sum_{\mathbf{v}_j \in B} n_j$, $I = \sum_{j=1}^{281} m_j$ and $J = \sum_{j=1}^{281} n_j$. Gangnon and Clayton [12] considered the ALR test statistic $U_{\mathcal{B}}^{(2)}$ with

$$\mathcal{B} = \{C(\mathbf{v}_i, r_{ij}) : 0 \leq r_{ij} \leq 20, 1 \leq i \leq 281, 1 \leq j \leq 281\},$$

where $r_{ij} = \|\mathbf{v}_i - \mathbf{v}_j\|$. We plot in Figure 2 simulated values of $U_{\mathcal{B}}^{(2)}$ under the null hypothesis (2.1) with $p_0 = 5 \times 10^{-4} (\doteq I/J)$, against quantiles of the chi-square distribution with one degree of freedom and also against quantiles of a distribution function $G$ satisfying

$$(4.1) \qquad G(x) = 1 - \left(\frac{2e^{-x}}{\pi x}\right)^{1/2}$$

for $x \geq x_0$ with $x_0 \doteq 0.42$ satisfying $2e^{-x_0}/(\pi x_0) = 1$. The upper tail probabilities of $G$ are expressions often seen in large deviations saddlepoint approximations.

Since $P\{\chi_1^2 \geq x\} \leq 1 - G(x)$ for all $x \geq 0$ and $P\{\chi_1^2 \geq x\} \sim 1 - G(x)$ as $x \to \infty$, $p$-value estimates of $U_{\mathcal{B}}^{(k)}$ obtained by comparing against $G$ instead of the chi-square distribution are slightly more conservative for small $p$-values. From the qq-plots, we see that $G$ provides a good fit over a wider range of values but for small $p$-values, which are of primary interest, the $p$-value estimates are comparable.

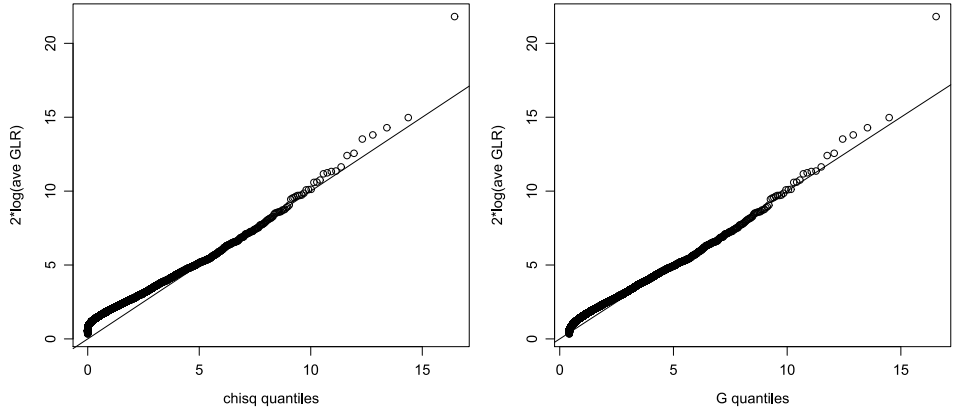

FIG. 2.  *Qq-plots of simulated values of $U_{\mathcal{B}}^{(2)}$ against the chi-square distribution with one degree of freedom (left) and the distribution $G$ (right).*

4.3. *Simulated datasets.*  The example in Section 4.2 is typical for application of cluster detection methodology. More than 20,000 circles were created from comparisons among the 281 census tracts. For larger number of census tracts, the number of circles can easily run into the millions. The computational burden is quite serious if say $L = 999$ or 9999 Monte Carlo simulation runs are used to evaluate $p$-values. Small $p$-values are of statistical interest, yet it is precisely for these cases that Monte Carlo methods are less reliable. If a person is looking at multiple regions, end-points or time-points, nominal $p$-values much smaller than 0.01 may be required for significance to be declared. For probability 0.05, $L = 999$ runs will give us relative error of about 0.15, while the corresponding relative error is about 0.3 for a probability 0.01. In Example 3 below, we compare the analytical chi-square and $G$ tail approximations [see (4.1)] of the ALR for two different arrangements of scanning sets. The key advantage of the analytical approximations lies in composite null situations for which the usual Monte Carlo methods may not work well (see Section 3).

EXAMPLE 3.  Let $\mathbf{v}_1, \ldots, \mathbf{v}_n$ be generated uniformly from the unit square $[0,1]^2$, and let

$$
\begin{aligned}
(4.2) \qquad \mathcal{B}_1 = \{ C(\mathbf{v}_i, r_{ij}) : 0 \le r_{ij} \le w_1, 1 \le i \le n, 1 \le j \le n \} \\
\text{where } r_{ij} = \| \mathbf{v}_i - \mathbf{v}_j \|.
\end{aligned}
$$

We shall abuse notation here and denote $\#\{i : \mathbf{v}_i \in C\}$ more simply by $\#C$.

For each $1 \le \ell \le L$ with $L$ large, generate independent standard normal random variables $Y_{1\ell}, \ldots, Y_{n\ell}$ and define

$$
Z_{C\ell} = \frac{\sum_{\mathbf{v}_i \in C}(Y_{i\ell} - \bar{Y}_\ell)}{\sqrt{(\#C)[1 - (\#C)/n]}} \qquad \text{where } \bar{Y}_\ell = n^{-1} \sum_{i=1}^{n} Y_{i\ell}.
$$

Let $U_{Z,\ell}^{(2)} = 2\log((\#\mathcal{B})^{-1}\sum_{C\in\mathcal{B}}e^{Z_{C\ell}^2/2})$. In Figure 3, we plot ordered values of $U_{Z,\ell}^{(2)}$ against quantiles of both the chi-square and $G$ distributions for $w_1 = 0.2$ and various values of $n$. Approximately six hours of computer time were taken up to generate the plot for $n = 1000$. The plots show the $G$ distribution to be more suitable for estimating moderately small $p$-values. For smaller $p$-values, the chi-square and $G$ distributions give similar approximations. Similar plots are obtained when experimenting with $w_1 = 0.3$.

**5. Discussion.** The New York leukaemia dataset in Section 4.2 is a typical dataset in which the locations are concentrated on a number of geographical centers instead of spreading over a domain $D$, and strictly speaking, the positive density assumption [see (A2)] does not hold. However, the purpose of the assumption is to ensure that the number of subjects in each scanning set goes to infinity at a fast enough rate, and as this is satisfied in this situation, the chi-square approximation is still valid. Similarly, the restriction that the class of sets in (A3) has to be all of the same shape can be relaxed in these types of datasets. The relaxation allows us to deal with the detection of irregular shaped clusters considered in, for example, [25, 30]. The condition that $\mathcal{B}$ be dependent only on the locations $\mathbf{x}_i$ and not on the responses $X_i$ is, however, necessary for the chi-square approximation to hold.

The qq-plots in Figures 2 and 3 show that the $p$-value approximations are inaccurate for small thresholds. This is consistent with the conditions of Theorem 1, which says that moderate or larger values of the threshold are needed for the $p$-value approximations to be accurate. This is not a problem because when large $p$-values are encountered, it suffices to state that the $p$-value is larger than a specified significance level. For very large thresholds, the difference of the approximated and empirical quantiles is due to the inaccuracy of the Monte Carlo method. Though Theorem 1 is stated only in terms of approximating unconditional $p$-values, a rough calculation shows that the chi-square approximation on the four conditional $p$-values of the ALR test statistics in Table 3 has the accuracy of about 4000 simulation runs. The chi-square approximations are also within one standard error of the Monte Carlo $p$-values in Table 5.

The overfitting of nuisance parameters when using Monte Carlo methods for $p$-value estimation of the spatial scan statistic was mentioned by Neill, Moore and Cooper [24], and this phenomenon likely contributed to the conservative $p$-values seen in Examples 1 and 2. The authors provided convincing arguments for why quick detection of disease outbreaks is important and cited the need to perform time-consuming Monte Carlo or bootstrap replications to provide reliable $p$-values of the spatial scan statistic as one justification for developing alternative methodologies. In this paper, we stick to the method of detection cluster via GLR values (but taking averages instead of maximums) popularized by the SaTScan software and address its drawbacks by providing accurate and easy to compute $p$-value approximations. These tail probability approximations can be applied even when nuisance parameters are in the model, and they enhance the attractiveness of the GLR method by easing its use.

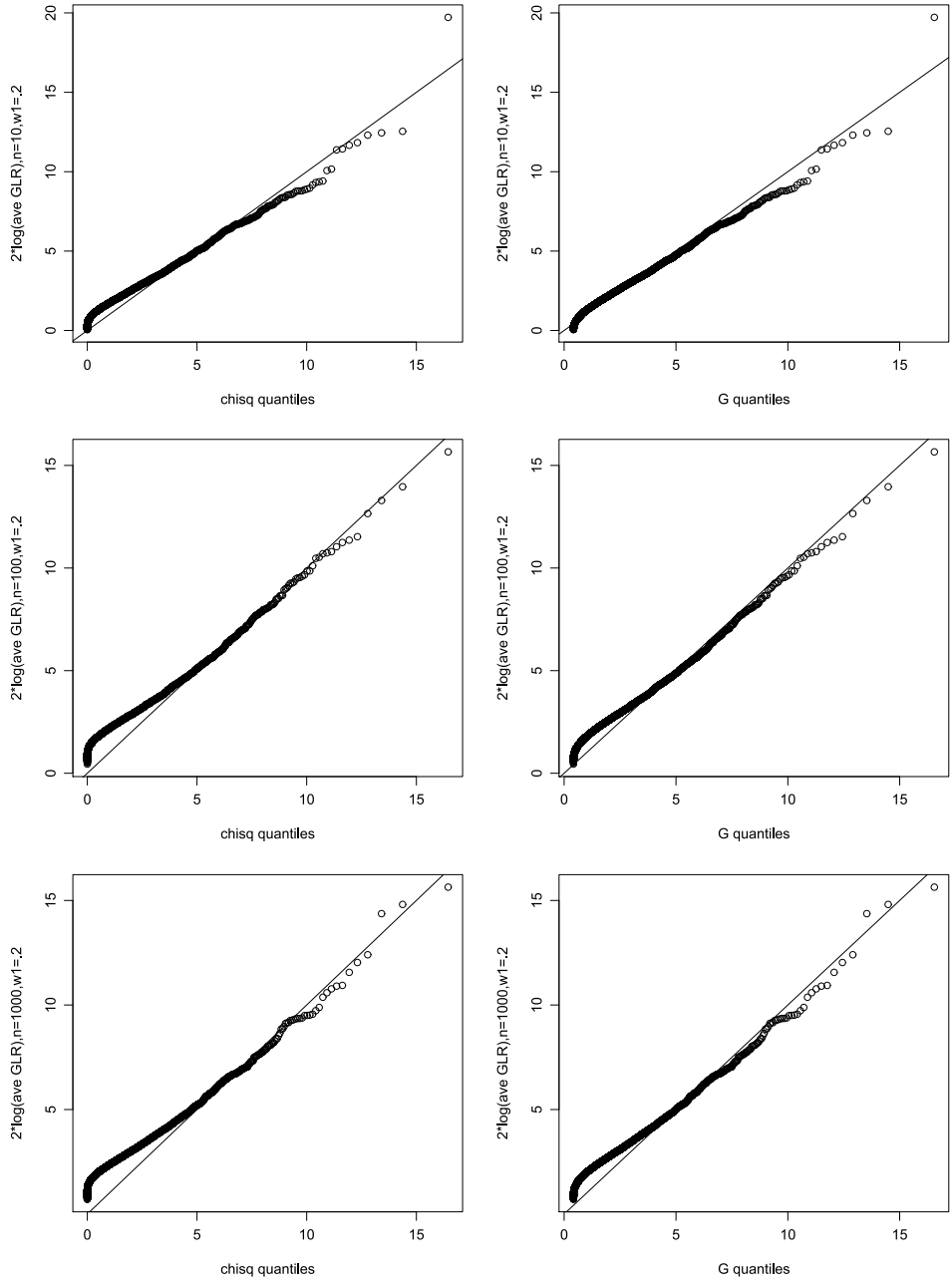

Fig. 3. *Qq-plots of $U_Z^{(2)}$ against the chi-square and G distributions for $\mathcal{B}_1$ with $n = 10$, 100, 1000 locations and maximum radius $w_1 = 0.2$.*

## 6. Proofs.

6.1. *Proof of Theorem 1.* Let $0 < \gamma_0 < p_0 < \gamma_1 < 1$. Then by large deviations,

$$(6.1) \qquad P_0\{\widehat{p}_0 \leq \gamma_0\} + P_0\{\widehat{p}_0 \geq \gamma_1\} = o(c^{-1/2}e^{-c/2}),$$

while by the law of large numbers, we may assume that

$$(6.2) \qquad \liminf_{J \to \infty} \inf_{B \in \mathcal{B}} (n_B/J) > 0.$$

For each $J\gamma_0 \leq I \leq J\gamma_1$, let $(p_B, \widetilde{p}_B)$ be the roots $(p, \widetilde{p})$ of

$$(6.3) \qquad \begin{aligned} n_B p + (J - n_B)\widetilde{p} &= I, \\ n_B\phi(p) + (J - n_B)\phi(\widetilde{p}) &= c/2 \qquad \text{with } p > \widehat{p}_0. \end{aligned}$$

Under (6.2), $(p_B, \widetilde{p}_B)$ exists and are unique for all $B \in \mathcal{B}$ when $J$ is large.

For given values of $\widehat{p}_0$ and $\mathbf{x}$, let $Q_B$ be a probability measure under which $X_1, \ldots, X_J$ are independent Bernoulli random variables satisfying

$$(6.4) \qquad Q_B\{X_i = 1 | \mathbf{x}_i \in B\} = p_B, \qquad Q_B\{X_i = 1 | \mathbf{x}_i \notin B\} = \widetilde{p}_B.$$

Let $\theta(p) = \log(p/\widehat{p}_0) - \log[(1-p)/(1-\widehat{p}_0)]$. Then by (2.3) and (6.3),

$$(6.5) \qquad \begin{aligned} \ell(B) &:= \log\left[\frac{dQ_B}{dP_{\widehat{p}_0}}(\mathcal{X})\right] \\ &= \sum_{\mathbf{x}_i \in B} \left\{\theta(p_B)X_i + \log\left(\frac{1-p_B}{1-\widehat{p}_0}\right)\right\} \\ &\quad + \sum_{\mathbf{x}_i \notin B} \left\{\theta(\widetilde{p}_B)X_i + \log\left(\frac{1-\widetilde{p}_B}{1-\widehat{p}_0}\right)\right\} \\ &= c/2 + \theta(p_B)\sum_{\mathbf{x}_i \in B}(X_i - p_B) + \theta(\widetilde{p}_B)\sum_{\mathbf{x}_i \notin B}(X_i - \widetilde{p}_B). \end{aligned}$$

The following supporting lemmas hold uniformly over $\gamma_0 \leq \widehat{p}_0 \leq \gamma_1$ under the conditions of Theorem 1. The proof of Lemma 1 is given in Appendix C, while the proof of Lemma 2(a) uses arguments in the proof of (6.11) which is also given in Appendix C. The proof of Lemma 2(b) is relatively straightforward and thus omitted.

LEMMA 1. *Assume* (6.2).
(a) $S^{(1)}(B) \geq \ell(B)$ *for all* $B \in \mathcal{B}$.
(b) *There exists* $\eta_c \to 0$ *as* $c \to \infty$ *such that*

$$|S^{(1)}(B) - \ell(B)| \leq \eta_c \qquad \text{whenever } |S^{(1)}(B) - c/2| \leq c^{1/3}.$$

LEMMA 2. (a) *Let*

$$V_B = \log\left\{\sum_{C \in \mathcal{B}} (dQ_C/dQ_B)(\mathcal{X})\right\} = \log\left(\sum_{C \in \mathcal{B}} e^{\ell(C) - \ell(B)}\right).$$

*Then whenever* $c \leq U_{\mathcal{B}}^{(1)} \leq c + c^{1/3}$,

$$\ell(B) + V_B - \log(\#\mathcal{B}) = \log\left((\#\mathcal{B})^{-1} \sum_{C \in \mathcal{B}} e^{\ell(C)}\right) = U_{\mathcal{B}}^{(1)}/2 + o(1).$$

(b) $Q_B\{\sum_{i=1}^J X_i = I | \mathbf{x}\} \sim P_{\hat{p}_0}\{\sum_{i=1}^J X_i = I | \mathbf{x}\}$ *uniformly over* $B \in \mathcal{B}$.

We shall now provide the key arguments in the proof of Theorem 1. Let $B_{\max}$ maximizes $\sum_{\mathbf{x}_i \in B}(X_i - p_B)$ over $B \in \mathcal{B}$, with an arbitrary ordering imposed on $\mathcal{B}$ to break ties. Under $Q_B$, conditioned on $c \leq U_{\mathcal{B}}^{(1)} \leq c + c^{1/3}$ and $B_{\max} = B$, $\ell(B)$ has an asymptotic density $(2\pi c)^{-1/2}$ on the interval $(\frac{c - c^{1/3}}{2}, \frac{c + c^{1/3}}{2})$, and is asymptotically independent of both $V_B$ and $\mathbf{I}_{\{B_{\max} = B\}}$. The random variable $V_B$ summarizes information on the local fluctuations of the GLR values for sets near $B$ when $B_{\max} = B$, and its value is determined chiefly by a small set of $X_i$ with $\mathbf{x}_i$ near the boundary of $B$, because under $Q_B$, $e^{\ell(C) - \ell(B)}$ is small for $C$, far from $B$. Similarly, $\mathbf{I}_{\{B_{\max} = B\}}$ is determined by the values of $X_i$ with $\mathbf{x}_i$ located near the boundary of $B$. The test statistic $\ell(B)$, on the other hand, is asymptotically $N(c/2, c)$ under $Q_B$ and is asymptotically independent of any small set of $X_i$. We thus obtain formally, for $\gamma_0 \leq \hat{p}_0 \leq \gamma_1$,

$$
\begin{aligned}
&P_{\hat{p}_0}\{U_{\mathcal{B}}^{(1)} \geq c | I, \mathbf{x}\} \\
&\quad = \sum_{B \in \mathcal{B}} P_{\hat{p}_0}\{U_{\mathcal{B}}^{(1)} \geq c, B_{\max} = B | I, \mathbf{x}\} \\
&\quad \sim \sum_{B \in \mathcal{B}} E_{Q(B)}(e^{-\ell(B)}\mathbf{I}_{\{U_{\mathcal{B}}^{(1)} \geq c, B_{\max} = B\}} | I, \mathbf{x}) \\
&\quad \sim \sum_{B \in \mathcal{B}} E_{Q(B)}(E[e^{-\ell(B)}\mathbf{I}_{\{\ell(B) \geq c/2 - V_B + \log(\#\mathcal{B}), B_{\max} = B\}} | V_B] | I, \mathbf{x}) \\
&\quad \sim \sum_{B \in \mathcal{B}} E_{Q(B)}\left(\mathbf{I}_{\{B_{\max} = B\}} \int_{c/2 - V_B + \log(\#\mathcal{B})}^{\infty} (2\pi c)^{-1/2} e^{-y}\, dy \Big| I, \mathbf{x}\right)
\end{aligned}
$$

(6.6)

$$= (2\pi c)^{-1/2} e^{-c/2} (\#\mathcal{B})^{-1} \sum_{B \in \mathcal{B}} E_{Q(B)}(e^{V_B} \mathbf{I}_{\{B_{\max} = B\}} | I, \mathbf{x})$$

$$= (2\pi c)^{-1/2} e^{-c/2} (\#\mathcal{B})^{-1} \sum_{B \in \mathcal{B}} \sum_{C \in \mathcal{B}} Q_C \{B_{\max} = B | I, \mathbf{x}\}.$$

We then switch the summation signs in the last line of (6.6) to show that $P_{\hat{p}_0}\{U_{\mathcal{B}}^{(1)} \geq c | I, \mathbf{x}\} \sim (2\pi c)^{-1/2} e^{-c/2}$, and (2.7) for $k = 1$ then follow from (6.1). By (6.5), $P_{\hat{p}_0}(A | \mathbf{x}) = E_{Q(B)}(e^{-\ell(B)} \mathbf{I}_A | \mathbf{x})$, where $A = \{U_{\mathcal{B}}^{(1)} \geq c, B_{\max} = B, \sum_{i=1}^J X_i = I\}$, and the relation between the first and second lines of (6.6) follows from Lemma 2(b). For additional details on (6.6), see Appendix C.

Since $S^{(2)}(B) = S^{(2)}(D \setminus B)$, and $\hat{p}_0$ lies between $m_B/n_B$ and $m_{D \setminus B}/n_{D \setminus B}$, it follows that

$$(6.7) \qquad e^{S^{(2)}(B)} = e^{S^{(1)}(B)} + e^{S^{(1)}(D \setminus B)} - 1.$$

Let $\widetilde{\mathcal{B}} = \{D \setminus B : B \in \mathcal{B}\}$ and $U_{\mathcal{B} \cup \widetilde{\mathcal{B}}}^{(1)} = 2 \log([2(\#\mathcal{B})]^{-1} \sum_{B \in \mathcal{B} \cup \widetilde{\mathcal{B}}} e^{S^{(1)}(B)})$. Then by the arguments leading to (2.7) for $k = 1$,

$$(6.8) \qquad P\{U_{\mathcal{B} \cup \widetilde{\mathcal{B}}}^{(1)} \geq c - 2 \log 2\} \sim [2\pi(c - 2\log 2)]^{-1/2} e^{-(c - 2\log 2)/2}$$

$$\sim [2/(\pi c)]^{1/2} e^{-c/2}.$$

By (6.7), $U_{\mathcal{B}}^{(2)} = 2 \log((\#\mathcal{B})^{-1} \sum_{B \in \mathcal{B}} e^{S^{(2)}(B)}) = U_{\mathcal{B} \cup \widetilde{\mathcal{B}}}^{(1)} + 2 \log 2 + o(1)$ when $U_{\mathcal{B}}^{(2)} \geq c$, and hence (2.7) for $k = 2$ follows from (6.8).

6.2. *Proof of Theorem 2.* Let $\Sigma = (\rho_{BC})_{B,C \in \mathcal{B}}$ be the covariance matrix of $Z = (Z_B)_{B \in \mathcal{B}}$, a multivariate normal with $EZ_B = 0$ and $\mathrm{Var}(Z_B) = 1$ for all $B \in \mathcal{B}$ under probability measure $P$. Hence $\rho_{BC} = \lambda_{B \cap C} - \lambda_B \lambda_C$. Fix $c > 0$, $B \in \mathcal{B}$ and let $Q_B$ be a probability measure under which $\omega(A) \sim N(\lambda_B^{-1/2}(1 - \lambda_B)^{1/2} \lambda_A c^{1/2}, \lambda_A)$ for $A \subset B$ and $\omega(A) \sim N(-\lambda_B^{1/2}(1 - \lambda_B)^{-1/2} \lambda_A c^{1/2}, \lambda_A)$ for $A \subset D \setminus B$, with $\omega(A), \omega(C)$ independent when $A$ and $C$ are disjoint sets. Under $Q_B$, $Z$ is multivariate normal with covariance matrix $\Sigma$ and $E_{Q(B)} Z_C = c^{1/2} \rho_{BC}$ for all $C \in \mathcal{B}$. Moreover,

$$(6.9) \qquad \ell(B) := \log\left[\frac{dQ_B}{dP}(Z)\right] = c^{1/2} Z_B - c/2.$$

We next use a linearization argument to justify the replacement of $Z_{B+}^2/2$ in the expression of $U_Z^{(1)}$ by $\ell(B)$. By convexity, $Z_{B+}^2/2 \geq \ell(B)$ for all $B \in \mathcal{B}$, with equality when $Z_{B+} = c^{1/2}$. By a Taylor expansion,

$$(6.10) \qquad \sup_{Z_B : |Z_{B+}^2 - c| \leq 2c^{1/3}} |\ell(B) - Z_{B+}^2/2| = O(c^{-1/3}) \qquad \text{as } c \to \infty.$$

Let $V_B = \log\{\sum_{C \in \mathcal{B}}(dQ_C/dQ_B)(Z)\} = \log(\sum_{C \in \mathcal{B}} e^{\ell(C) - \ell(B)})$. Then by (6.10), there exists $\zeta_c = O(c^{-1/3})$ such that whenever $c \leq U_Z^{(1)} \leq c + c^{1/3}$,

$$U_Z^{(1)}/2 - \zeta_c \leq \ell(B) + V_B - \log(\#\mathcal{B})\left[ = \log\left((\#\mathcal{B})^{-1}\sum_{C \in \mathcal{B}} e^{\ell(C)}\right)\right]$$

(6.11)
$$\leq U_Z^{(1)}/2,$$

(see Appendix C). We then apply the steps in (6.6), without the conditioning on $I$ and $\mathbf{x}$, to obtain the tail probabilities of $U_Z^{(1)}$. For extensions to the tail probabilities of $U_Z^{(2)}$, apply the arguments in the last paragraph of Section 6.1.

## APPENDIX A: THEOREM 3 AND ITS PROOF

THEOREM 3. *Let $S_{1c}, \ldots, S_{nc}$ be random variables and assume that there exists a constant $K > 0$ and random variables $Y_{kj}$ such that $P\{Y_{kj} = 0\} = 0$ for all $k \neq j$ and as $c \to \infty$,*

$$P\{S_{kc} \geq c + y\} \sim Kc^{-1/2}e^{-c-y}, \tag{A.1}$$

*while conditioned on $S_{kc} \geq c + y$,*

$$(S_{kc} - S_{1c}, \ldots, S_{kc} - S_{nc}) \Rightarrow (Y_{k1}, \ldots, Y_{kn}) \tag{A.2}$$

*for all $1 \leq k \leq n$ and $y \in \mathbf{R}$. Then*

$$n^{-1}\sum_{k=1}^{n} E\left[\left(\sum_{j=1}^{n} e^{-Y_{kj}}\right)\mathbf{I}_{\{Y_{kj} \geq 0 \text{ for all } j\}}\right] = 1 \tag{A.3}$$

*and*

$$P\left\{n^{-1}\sum_{j=1}^{n} e^{S_{jc}} \geq e^c\right\} \sim Kc^{-1/2}e^{-c} \qquad \text{as } c \to \infty. \tag{A.4}$$

PROOF. Let $M_c = \sup_{1 \leq k \leq n} S_{kc}$. For a given $\varepsilon > 0$, let $0 = y_1 < \cdots < y_m$ be such that $P\{Y_{kj} = y_r\} = 0$ for all $1 \leq r \leq m$, $k \neq j$ and $\sup_{1 \leq r \leq m}(e^{-y_r} - e^{-y_{r+1}}) \leq \varepsilon$, where $y_{m+1} = \infty$. Then by (A.1) and (A.2), for all $k \neq j$,

$$P\{S_{jc} \geq c, M_c = S_{kc}\}$$

$$= \sum_{r=1}^{m} P\{S_{jc} \geq c, M_c = S_{kc}, y_r \leq S_{kc} - S_{jc} < y_{r+1}\}$$

$$\leq \sum_{r=1}^{m} P\{S_{kc} \geq c + y_r, M_c = S_{kc}, y_r \leq S_{kc} - S_{jc} < y_{r+1}\} \tag{A.5}$$

$$\sim K c^{-1/2} e^{-c} \sum_{r=1}^{m} e^{-y_r} P\{Y_{ki} \geq 0 \text{ for all } i, y_r \leq Y_{kj} < y_{r+1}\}$$

$$\leq K c^{-1/2} e^{-c} E[(e^{-Y_{kj}} + \varepsilon) \mathbf{I}_{\{Y_{ki} \geq 0 \text{ for all } i\}}].$$

Similarly,

$$P\{S_{jc} \geq c, M_c = S_{kc}\}$$

(A.6)
$$\geq (K + o(1)) c^{-1/2} e^{-c} \sum_{r=1}^{m} e^{-y_{r+1}} P\{Y_{ki} \geq 0 \text{ for all } i, y_r \leq Y_{kj} < y_{r+1}\}$$

$$\geq (K + o(1)) c^{-1/2} e^{-c} E[(e^{-Y_{kj}} - \varepsilon) \mathbf{I}_{\{Y_{ki} \geq 0 \text{ for all } i\}}].$$

By selecting $\varepsilon$ arbitrarily small, it follows from (A.5) and (A.6) that

(A.7) $$P\{S_{jc} \geq c, M_c = S_{kc}\} \sim K c^{-1/2} e^{-c} E(e^{-Y_{kj}} \mathbf{I}_{\{Y_{ki} \geq 0 \text{ for all } i\}}).$$

The asymptotic relation (A.7) also holds for $k = j$ by applying (A.1) and (A.2) for $y = 0$, noting that $Y_{jj}$ is a zero-valued random variable for all $j$. We then add up (A.7) over $1 \leq j \leq n$, $1 \leq k \leq n$ and compare against the asymptotic relation $\sum_{j=1}^{n} P\{S_{jc} \geq c\} \sim K n c^{-1/2} e^{-c}$, which follows from (A.1), to obtain (A.3).

Since $\log(n^{-1} \sum_{j=1}^{n} e^{S_{jc} - S_{kc}}) \leq 0$ when $M_c = S_{kc}$, by (A.1) and (A.2),

$$P\left\{ n^{-1} \sum_{j=1}^{n} e^{S_{jc}} \geq e^c, M_c = S_{kc} \right\}$$

(A.8)
$$= P\left\{ S_{kc} \geq c - \log\left( n^{-1} \sum_{j=1}^{n} e^{S_{jc} - S_{kc}} \right), M_c = S_{kc} \right\}$$

$$\sim K c^{-1/2} e^{-c} E\left[ \left( n^{-1} \sum_{j=1}^{n} e^{-Y_{kj}} \right) \mathbf{I}_{\{Y_{kj} \geq 0 \text{ for all } j\}} \right]$$

and (A.4) follows from (A.3) by adding (A.8) over $1 \leq k \leq n$. To show the last relation in (A.8), we use a discretization argument described earlier. Given any $\varepsilon > 0$, select $0 = y_1 < \cdots < y_m$ such that $P\{-\log(n^{-1} \sum_{j=1}^{n} e^{-Y_{kj}}) = y_r, Y_{kj} \geq 0 \text{ for all } j\} = 0$ for all $1 \leq r \leq m$, $1 \leq k \leq n$ and $\sup_{1 \leq r \leq m}(e^{-y_r} - e^{-y_{r+1}}) \leq \varepsilon$, with $y_{m+1} = \infty$. We then express asymptotic upper and lower bounds of

$$P\left\{ S_{kc} \geq c - \log\left( n^{-1} \sum_{j=1}^{n} e^{S_{jc} - S_{kc}} \right), \right.$$

$$\left. M_c = S_{kc}, y_r \leq -\log\left( n^{-1} \sum_{j=1}^{n} e^{S_{jc} - S_{kc}} \right) < y_{r+1} \right\},$$

in terms of $\varepsilon$ and expectations involving $Y_{ki}$ before letting $\varepsilon \to 0$. The details are omitted. $\square$

## APPENDIX B: ASYMPTOTIC EXPANSIONS OF THE LOG LIKELIHOOD FUNCTION

For a fixed $B \in \mathcal{B}$, let $z_i(\theta, \beta) = \beta' \mathbf{u}_i + \theta \mathbf{I}_{\{\mathbf{x}_i \in B\}}$ and let

$$\ell_i(\theta, \beta) = -X_i \log(1 + e^{-z_i(\theta,\beta)}) - (1 - X_i) \log(1 + e^{z_i(\theta,\beta)})$$

be the log likelihood function corresponding to the $i$th subject. Since $\partial \ell_i / \partial z_i = X_i - p_i$, where $p_i = (1 + e^{-z_i})^{-1}$, evaluated at some parameter $\beta$ and $\theta = 0$, it follows that

$$\frac{d\ell_i}{d\theta} = (X_i - p_i)\mathbf{I}_{\{\mathbf{x}_i \in B\}} \quad \text{and} \quad \frac{d\ell_i}{d\beta_k} = u_{ik}(X_i - p_i) \qquad \text{for } 1 \le k \le r.$$

To motivate the form of the limiting distribution of $S^{(k)}(B)$, we use a weighted Gram–Schmidt orthogonalization procedure, rather than matrix notation, to describe the first-order quadratic term in a Taylor expansion of $S^{(2)}(B)$. Let $w_i = p_i(1 - p_i)$ and let weighted dot product $(\mathbf{a} \cdot \mathbf{b})_w = \sum_{i=1}^{J} a_i b_i w_i$ and norm $\|\mathbf{a}\|_w = (\mathbf{a} \cdot \mathbf{a})_w$. Define recursively $\widetilde{\mathbf{u}}_1 = \mathbf{u}_1$ and $\widetilde{\mathbf{u}}_k = \mathbf{u}_k - \sum_{s=1}^{k-1} a_{ks} \widetilde{\mathbf{u}}_s$, where $a_{ks} = (\mathbf{u}_k \cdot \widetilde{\mathbf{u}}_s)_w / \|\widetilde{\mathbf{u}}_s\|_w^2$, for $2 \le k \le r$. Then $(\widetilde{\mathbf{u}}_k \cdot \widetilde{\mathbf{u}}_s)_w = 0$ for all $k \ne s$. Let $\widetilde{\mathbf{u}}_k = (\widetilde{u}_{1k}, \ldots, \widetilde{u}_{Jk})$. Under sufficient regularity conditions, $S_B^{(2)}$ is equal to, up to a $o(1)$ term,

$$(2v_B^2)^{-1} \left\{ \sum_{i=1}^{J} \left( \mathbf{I}_{\{\mathbf{x}_i \in B\}} - \sum_{k=1}^{r} \frac{(\alpha_B \cdot \widetilde{\mathbf{u}}_k)_w \widetilde{u}_{ik}}{\|\widetilde{\mathbf{u}}_k\|_w^2} \right) (X_i - p_i) \right\}^2,$$

(B.1)

$$\text{where } \alpha_B = (\mathbf{I}_{\{\mathbf{x}_1 \in B\}}, \ldots, \mathbf{I}_{\{\mathbf{x}_J \in B\}})'$$

$$\text{and } v_B^2 = \sum_{i=1}^{J} w_i \left( \mathbf{I}_{\{\mathbf{x}_i \in B\}} - \sum_{k=1}^{r} \frac{(\alpha_B \cdot \widetilde{\mathbf{u}}_k)_w \widetilde{u}_{ik}}{\|\widetilde{\mathbf{u}}_k\|_w^2} \right)^2.$$

We will next consider a characterization of the limiting distributions of $S^{(2)}(B)$, $B \in \mathcal{B}$. Let $\eta(\mathbf{t}) = \lambda(\mathbf{t}) E(w_1 | \mathbf{x}_1 = \mathbf{t}) / E(w_1)$ and $g_k(\mathbf{t}) = E(u_{1k} w_1 | \mathbf{x}_1 = \mathbf{t}) / E(w_1 | \mathbf{x}_1 = \mathbf{t})$ for $1 \le k \le r$ and assume that they are positive and continuous on $D$. Let $\widetilde{g}_1 = g_1(=1)$ and define recursively for $k \ge 2$,

$$\widetilde{g}_k(\mathbf{t}) = g_k(\mathbf{t}) - \sum_{s=1}^{k-1} \mu_{ks} \widetilde{g}_s(\mathbf{t}),$$

where

$$\mu_{ks} = b_s^{-1} \int_D g_k(\mathbf{t}) \widetilde{g}_s(\mathbf{t}) \eta(\mathbf{t}) \, d\mathbf{t} \quad \text{and} \quad b_s = \int_D \widetilde{g}_s^2(\mathbf{t}) \eta(\mathbf{t}) \, d\mathbf{t}.$$

Then $\int_D \tilde{g}_k(\mathbf{t})\tilde{g}_s(\mathbf{t})\eta(\mathbf{t})\,d\mathbf{t} = 0$ for all $k \neq s$. Let $\omega$ be Gaussian white noise on $D$ with $\omega(B) \sim N(0, \eta_B)$, where $\eta_B = \int_B \eta(\mathbf{t})\,d\mathbf{t}$. Then $S^{(2)}(B)$ converges weakly to $Z_B^2/2$, where

$$Z_B = v_B^{-1}\left[\omega(B) - \sum_{k=1}^r \left(b_k^{-1}\int_B \tilde{g}_k(\mathbf{t})\eta(\mathbf{t})\,d\mathbf{t}\right)\omega_k(D)\right]$$

with $\omega_k(D) = \int_D \tilde{g}_k(\mathbf{t})\omega(d\mathbf{t})$ and $v_B$ is a normalizing constant to ensure $\mathrm{Var}(Z_B) = 1$.

The justification behind (3.6) requires arguments used in the proof of Theorem 2. For a given $c > 0$ and $B \in \mathcal{B}$, let $Q_B$ be a probability measure such that $\omega(A) \sim N(\int_A \mu(\mathbf{t})\eta(\mathbf{t})\,d\mathbf{t}, \eta_A)$, where $\mu(\mathbf{t}) = c^{1/2}v_B^{-1}[\mathbf{I}_{\{\mathbf{t}\in B\}} - \sum_{k=1}^r \gamma_{kB}\tilde{g}_k(\mathbf{t})]$ and $\gamma_{kB} = b_k^{-1}\int_B \tilde{g}_k(\mathbf{t})\eta(\mathbf{t})\,d\mathbf{t}$. Moreover, under $Q_B$, $\omega(A)$ and $\omega(C)$ are independent whenever $A$ and $C$ are disjoint. By Girsanov's theorem (see Chapter 3.5 of [14]),

$$\ell(B) := \log\left[\frac{dQ_B}{dP}(\omega)\right] = \int_D \mu(\mathbf{t})\omega(d\mathbf{t}) - \frac{1}{2}\int_D \mu^2(\mathbf{t})\eta(\mathbf{t})\,d\mathbf{t} = c^{1/2}Z_B - c/2.$$

We can then proceed, as in the proof of Theorem 2 in Section 6.2, by analyzing the behavior of $U_Z^{(k)}$ under $Q_B$ and using a linearization argument to estimate $Z_B^2/2$ by $\ell(B)$ when $Z_B$ is close to $c^{1/2}$. The details are omitted.

## APPENDIX C: PROOFS OF LEMMA 1, (6.11) AND (6.6)

PROOF OF LEMMA 1. Let $\alpha_B = n_B/J$ and $f(p) = \alpha_B\phi(p) + (1-\alpha_B) \times \phi((\hat{p}_0 - \alpha_B p)/(1-\alpha_B))$. The tangent of $f$ at $p = p_B$ is

$$g(p) := \alpha_B\{\theta(p_B)p + \log[(1-p_B)/(1-\hat{p}_0)]\}$$
$$+ (1-\alpha_B)\{\theta(\tilde{p}_B)(\hat{p}_0 - \alpha_B p)/(1-\alpha_B) + \log[(1-\tilde{p}_B)/(1-\hat{p}_0)]\}.$$

Since $S^{(1)}(B) = Jf(m_B/n_B)$ and $\ell(B) = Jg(m_B/n_B)$, Lemma 1(a) follows from the convexity of $f$.

Next, let $K = [p_B, p]$ if $p_B \leq p$ and $K = [p, p_B]$ if $p_B > p$. Since $f(p_B) = g(p_B)$ and $g$ is linear,

$$|f(p) - g(p)| \leq \left[\sup_{q\in K} f''(q)\right](p - p_B)^2/2,$$

(C.1)

$$|f(p) - f(p_B)| \geq \left[\inf_{q\in K} f'(q)\right]|p - p_B|.$$

Select $p = m_B/n_B$. Since $f'(\hat{p}_0) = f'(\hat{p}_0) = 0$ and $f(p_B) = c(2J)^{-1} = o(1)$, it follows that $f'(p_B)$ is of order $(c/J)^{1/2}$. If $J|f(p) - f(p_B)|(=|S(B) - c/2|) \leq c^{1/3}$, then by the first inequality of (C.1), $|p - p_B| = O(c^{-1/6}J^{-1/2})$. Since $|S^{(1)}(B) - \ell(B)| = J|f(p) - g(p)|$ and $f''$ is order 1 in $K$, Lemma 1(b) follows from second inequality of (C.1). $\square$

PROOF OF (6.11). The upper bound follows from $Z_{C+}^2/2 \geq \ell(C)$ for all $C$. Next, observe that the constraints $U_Z^{(1)} \leq c + c^{1/3}$ and $\log(\#\mathcal{B}) = o(c^{1/3})$ together imply that $\sup_{C \in \mathcal{B}} Z_{C+}^2 \leq c + 2c^{1/3}$ for all large $c$. Since

$$(\#\mathcal{B})^{-1} \sum_{C \in \mathcal{B}} (e^{Z_{C+}^2/2} \mathbf{I}_{\{Z_{C+}^2 < c - 2c^{1/3}\}}) = o(e^{c/2}),$$

the lower bound follows from applying (6.10) on

$$(\#\mathcal{B})^{-1} \sum_{C \in \mathcal{B}} (e^{Z_{C+}^2/2} \mathbf{I}_{\{|Z_{C+}^2 - c| \leq 2c^{1/3}\}}). \qquad \square$$

PROOF OF (6.6). Let $\partial B$ denote the boundary of $B$ and let

$$\partial_\varepsilon B = \{\mathbf{t} \in D : \|\mathbf{t} - \mathbf{u}\| \leq \varepsilon \text{ for some } \mathbf{u} \in \partial B\} \qquad \text{with } \varepsilon = c^{-3/5}.$$

Let $B \triangle C = (B \setminus C) \cup (C \setminus B)$ and let $\mathcal{B}_1 = \{C \in \mathcal{B} : B \triangle C \subset \partial_\varepsilon B\}$, the class of $C \in \mathcal{B}$ that are "close" to $B$. In Lemma 3(a) below, we show that $\ell(C) - \ell(B)$ can be approximated by

$$
\begin{aligned}
h_B(C) := &\sum_{\mathbf{x}_i \in C \setminus B} [\theta(p_C)(X_i - p_B) - \theta(\widetilde{p}_C)(X_i - \widetilde{p}_B)] \\
&- \sum_{\mathbf{x}_i \in B \setminus C} [\theta(p_C)(X_i - p_B) - \theta(\widetilde{p}_C)(X_i - \widetilde{p}_B)]
\end{aligned}
\tag{C.2}
$$

for all $C \in \mathcal{B}_1$. We show in Lemma 3(b) that $\sum_{C \notin \mathcal{B}_1} e^{\ell(C) - \ell(B)}$ is asymptotically negligible. Hence $V_B = \log(\sum_{C \in \mathcal{B}_1} e^{h_B(C)}) + o(1)$. But $\log(\sum_{C \in \mathcal{B}_1} e^{h_B(C)})$ depends only on $\mathcal{X}_B = \{(\mathbf{x}_i, X_i) : \mathbf{x}_i \in \partial_\varepsilon B\}$ and because $\varepsilon = o(c^{-1/2})$, $\ell(B)$ is asymptotically independent of $\mathcal{X}_B$. Let

$$\Gamma_{\beta_1, \beta_2}(B) = \{\mathbf{x} : \#(\partial_\varepsilon B) \leq \beta_1 J\varepsilon \text{ and } \#(B \triangle C) \geq \beta_2 J\varepsilon \text{ for all } C \in \mathcal{B} \setminus \mathcal{B}_1\}. \qquad \square$$

LEMMA 3. *There exists $\beta_1 > 0$ large enough and $\beta_2 > 0$ small enough such that*

$$1 - P(\Gamma_{\beta_1, \beta_2}(B)) = o(e^{-c^{1/3}}) \qquad \text{uniformly over } B \in \mathcal{B}. \tag{C.3}$$

*Moreover, for fixed $\beta_1 > 0$ and $\beta_2 > 0$, the following holds uniformly over $\mathbf{x} \in \Gamma_{\beta_1, \beta_2}(B)$.*

(a) *If $c/2 \leq \ell(B) \leq c/2 + c^{1/3}$ and $\sum_{i=1}^J X_i = I$, then*

$$\ell(C) - \ell(B) = h_B(C) + o(1) \qquad \text{uniformly over } C \in \mathcal{B}_1. \tag{C.4}$$

(b) *Let $\Lambda(C) = \{\ell(B) \geq (c/2) \text{ and } \ell(C) \geq (c/2) - c^{1/3}\}$. Then*

$$P\{\Lambda(C) | I, \mathbf{x}\} = o(e^{-c/2 - c^{1/3}}) \qquad \text{uniformly over } C \in \mathcal{B} \setminus \mathcal{B}_1. \tag{C.5}$$

PROOF.  Let $\sigma^{d-1}(\cdot)$ denote a $(d-1)$-dimensional volume element. By (A2), $\#(\partial_\varepsilon B) \sim \mathrm{Bin}(J, q)$, where $q \sim \varepsilon\zeta_B$ and $\zeta_B = 2\int_{\partial B} \lambda(\mathbf{t})\sigma^{d-1}(d\mathbf{t})$. Then

(C.6) $\quad P\{\#\partial_\varepsilon B \geq 3\zeta_B \varepsilon J\} = o(e^{-\zeta_B \varepsilon J \eta_0}) \qquad$ for some $\eta_0 > 0$.

Similarly, there exists $\kappa_B > 0$ such that for all $C \in \mathcal{B} \setminus \mathcal{B}_1$ and $J$ large, $\#(B \triangle C) \sim \mathrm{Bin}(J, q_C)$, where $q_C \geq \varepsilon\kappa_B$. Hence

(C.7) $\quad P\{\#(B \triangle C) \leq \kappa_B \varepsilon J/3\} = o(e^{-\kappa_B \varepsilon J \eta_1}) \qquad$ for some $\eta_1 > 0$.

Since $\varepsilon J/c^{1/3} \to \infty$ and $\log(\#\mathcal{B}) = o(c^{1/3})$, (C.3) follows from (C.6) and (C.7).

(a) By (6.3) and (6.5),

$$\ell(B) = n_B \phi(p_B) + (J - n_B)\phi(\widetilde{p}_B)$$
$$+ \theta(p_B) \sum_{\mathbf{x}_i \in B} (X_i - p_B) + \theta(\widetilde{p}_B) \sum_{\mathbf{x}_i \notin B} (X_i - \widetilde{p}_B),$$

$$\ell(C) = n_C\{p_B \log(p_C/\widehat{p}_0) + (1 - p_B)\log[(1 - p_C)/(1 - \widehat{p}_0)]\}$$
$$+ (J - n_C)\{\widetilde{p}_B \log(\widetilde{p}_C/\widehat{p}_0) + (1 - \widetilde{p}_B)\log[(1 - \widetilde{p}_C)/(1 - \widehat{p}_0)]\}$$
$$+ \theta(p_C) \sum_{\mathbf{x}_i \in C} (X_i - p_B) - \theta(\widetilde{p}_C) \sum_{\mathbf{x}_i \notin C} (X_i - \widetilde{p}_B).$$

If $c/2 \leq \ell(B) \leq c/2 + c^{1/3}$, then by (C.2),

$$\ell(C) - h_B(C) = n_C\{p_B \log(p_C/\widehat{p}_0) + (1 - p_B)\log[(1 - p_C)/(1 - \widehat{p}_0)]\}$$
$$+ (J - n_C)\{\widetilde{p}_B \log(\widetilde{p}_C/\widehat{p}_0)$$
$$+ (1 - \widetilde{p}_B)\log[(1 - \widetilde{p}_C)/(1 - \widehat{p}_0)]\}$$
$$+ \theta(p_C) \sum_{\mathbf{x}_i \in B} (X_i - p_B) + \theta(\widetilde{p}_C) \sum_{\mathbf{x}_i \notin B} (X_i - \widetilde{p}_B)$$
$$= \ell(B) + O(J(p_C - p_B)^2) + O(c^{1/3}|p_C - p_B|).$$

Under $\Gamma_{\beta_1, \beta_2}(B)$, if $C \in \mathcal{B}_1$, then by (6.3), $|p_C - p_B| = O(c^{1/2}|n_B^{-1/2} - n_C^{-1/2}|) = O(c^{1/2}|n_B - n_C|J^{-3/2}) = O(c^{1/2}J^{-1/2}\varepsilon)$ and we conclude (C.4).

(b) For given $\widehat{p}_0$ and $\mathbf{x}$ generated according to (A2), let $Q_{B,C}$ be a probability measure under which $X_1, \ldots, X_J$ are independent Bernoulli random variables satisfying

$$Q_{B,C}\{X_i = 1\} = (p_B \mathbf{I}_{\{\mathbf{x}_i \in B\}} + \widetilde{p}_B \mathbf{I}_{\{\mathbf{x}_i \notin B\}} + p_C \mathbf{I}_{\{\mathbf{x}_i \in C\}} + \widetilde{p}_C \mathbf{I}_{\{\mathbf{x}_i \notin C\}})/2.$$

By the AM $\geq$ GM inequality and (6.4),

(C.8) $\quad Q_{B,C}\{X_i = a\} \geq (Q_B\{X_i = a\}Q_C\{X_i = a\})^{1/2}, \qquad a = 0, 1.$

By (6.3) and the identity $(x+y)/2 - (x^{1/2}y^{1/2}) = (x^{1/2} - y^{1/2})^2/2$, there exists $\gamma > 0$ such that whenever $\mathbf{x}_i \in B \triangle C$,

(C.9) $\quad Q_{B,C}\{X_i = a\} \geq e^{c\gamma/J}(Q_B\{X_i = a\}Q_C\{X_i = a\})^{1/2}, \qquad a = 0, 1.$

Let $C \in \mathcal{B} \setminus \mathcal{B}_1$. By (C.8), (C.9) and the relation $P_{p_0}\{\sum_{i=1}^J X_i = I\} \sim Q_{B,C}\{\sum_{i=1}^J X_i = I\}$,

$$P\{\Lambda(C) | I, \mathbf{x}\} \leq (1+o(1))E_{Q(B,C)}(e^{-[\ell(B)+\ell(C)]/2 - (c\gamma/J)(\#B\triangle C)}\mathbf{I}_{\Lambda(C)}|I, \mathbf{x})$$

and (C.5) holds because under $\Gamma_{\beta_1,\beta_2}(B)$, $\#(B \triangle C) \geq \beta_2 J\varepsilon$ for $C \in \mathcal{B} \setminus \mathcal{B}_1$. $\square$

To show that the second and fifth lines of (6.6) are asymptotically equivalent, assume without loss of generality that $\mathbf{x} \in [\bigcap_{B \in \mathcal{B}} \Gamma_{\beta_1,\beta_2}(B)]$ for $\beta_1 > 0$ and $\beta_2 > 0$ satisfying (C.3). By (6.3), if $\sum_{\mathbf{x}_i \in B}(X_i - p_B) < 0$ for all $B \in \mathcal{B}$, then $S^{(1)}(B) < c/2$ for all $B \in \mathcal{B}$ and hence $U_{\mathcal{B}} < c$. The condition $U_{\mathcal{B}} \geq c$ thus ensures that $\sum_{\mathbf{x}_i \in B_{\max}}(X_i - p_{B_{\max}}) \geq 0$, and, consequently, $\ell(B_{\max}) \geq c/2$ [see (6.5)].

Let $\Omega_B = \{\ell(B) \geq c/2$ and $\ell(C) < (c/2) - c^{1/3}$ for all $C \in \mathcal{B} \setminus \mathcal{B}_1\}$ and $W_{\mathcal{B}} = \log((\#\mathcal{B})^{-1} \sum_{C \in \mathcal{B}} e^{\ell(C)})$. By Lemma 2(a) and (C.5), there exists $c' = c + o(1)$ such that

$$\begin{aligned}
&E_{Q(B)}(e^{-\ell(B)}\mathbf{I}_{\{U_{\mathcal{B}} \geq c, B_{\max}=B\}}|I, \mathbf{x}) \\
&\text{(C.10)} \quad \geq E_{Q(B)}(e^{-\ell(B)}\mathbf{I}_{\{c'+c^{1/3} \geq 2W_{\mathcal{B}} \geq c', B_{\max}=B\}}|\Omega_B, I, \mathbf{x}) \\
&\qquad + o((\#\mathcal{B})^{-1}c^{-1/2}e^{-c/2}).
\end{aligned}$$

Let $\mathcal{X}_B = \{(\mathbf{x}_i, X_i) : \mathbf{x}_i \in \partial_\varepsilon B\}$ and assume $\Omega_B$. Then $\mathbf{I}_{\{B_{\max}=B\}}$ is determined on knowing $\mathcal{X}_B$, and, in addition, by (C.4), $V_B = \log(\sum_{C \in \mathcal{B}_1} e^{h_B(C)}) + o(1)$. Moreover, by (C.2), $\log(\sum_{C \in \mathcal{B}_1} e^{h_B(C)})$ is determined on knowing $\mathcal{X}_B$ and $\mathbf{x}$. Hence there exists $c_* = c + o(1)$ such that

$$\begin{aligned}
&E_{Q(B)}(e^{-\ell(B)}\mathbf{I}_{\{c'+c^{1/3} \geq 2W_{\mathcal{B}} \geq c', B_{\max}=B\}}|\Omega_B, I, \mathbf{x}, \mathcal{X}_B) \\
&\text{(C.11)} \quad \geq (1+o(1))\mathbf{I}_{\{B_{\max}=B\}}(\#\mathcal{B})^{-1}e^{V_B} \\
&\qquad \times E_{Q(B)}(e^{-W_{\mathcal{B}}}\mathbf{I}_{\{c_*/2+c^{1/3} \geq W_{\mathcal{B}} \geq c_*/2\}}|\Omega_B, I, \mathbf{x}, \mathcal{X}_B).
\end{aligned}$$

It follows from a local limit theorem that under $Q_B$, $W_{\mathcal{B}}$ conditioned on $\Omega_B$ has an asymptotic density of $(2\pi c)^{-1/2}$ uniformly over $[c_*/2, c_*/2 + c^{1/3}]$. Replace this asymptotic density into (C.11), take expectation over $\mathcal{X}_B$ and substitute the remaining expression into (C.10) to obtain

$$\begin{aligned}
&E_{Q(B)}(e^{-\ell(B)}\mathbf{I}_{\{U_{\mathcal{B}} \geq c, B_{\max}=B\}}|I, \mathbf{x}) \\
&\text{(C.12)} \quad \geq (1+o(1))(2\pi c)^{-1/2}e^{-c/2} \\
&\qquad \times (\#\mathcal{B})^{-1}E_{Q(B)}(e^{V_B}\mathbf{I}_{\{B_{\max}=B\}}|I, \mathbf{x}) + o((\#\mathcal{B})^{-1}c^{-1/2}e^{-c/2}).
\end{aligned}$$

Similarly, $E_{Q(B)}(e^{-\ell(B)}\mathbf{I}_{\{U_B \geq c+c^{1/3}, B_{\max}=B\}}|I, \mathbf{x}) = o((\#\mathcal{B})^{-1}c^{-1/2}e^{-c/2})$ and (C.12) with the inequality reversed can be obtained. Hence the second and fifth lines of (6.6) are asymptotically equivalent.

**Acknowledgments.** The author would like to thank Professor P. J. Diggle for the use of the laryngeal-lung cancer datasets. He would also like to thank the reviewers for their comments which led to improvements of the manuscript.

Department of Statistics
and Applied Probability
National University of Singapore
6 Science Drive 2
Singapore 117546
E-mail: stachp@nus.edu.sg